\documentclass[11pt]{article}

\usepackage[margin=1in]{geometry}
\usepackage{amsmath,amssymb,amsthm,mathtools}
\usepackage{hyperref}
\usepackage{enumitem}
\usepackage{tikz}

\newtheorem{theorem}{Theorem}[section]
\newtheorem{lemma}[theorem]{Lemma}
\newtheorem{proposition}[theorem]{Proposition}
\newtheorem{corollary}[theorem]{Corollary}
\theoremstyle{remark}
\newtheorem{remark}[theorem]{Remark}
\newtheorem{problem}{Problem}

\newcommand{\codeg}{\operatorname{codeg}}

\newcommand{\eps}{\varepsilon}
\newcommand{\TK}{\mathrm{TK}}
\newcommand{\TKshort}{\mathrm{TK}^{(\leq 2)}}

\title{Topological Clique Exclusion in Positively Curved Graphs}
\author{Linyuan Lu\thanks{University of South Carolina, Columbia, SC 29208,
({\tt lu@math.sc.edu}).
This author was supported in part by NSF grants DMS 2038080 and DMS 2601952.
}}
\date{\today}

\begin{document}
\maketitle

\begin{abstract}
We study locally finite graphs of positive Lin--Lu--Yau curvature and minimum
degree at least $t-1$ under two exclusions: a topological $K_t$, denoted
$\TK_t$, and its locally topological version $\TKshort_t$, in which every
replacement path has length at most two. The ordinary exclusion gives an
$O(t^{7})$ maximum-degree bound by combining the classical density theorem
for topological cliques with a curvature-driven rainbow-clique argument.
For $K_t$-minor-free graphs, a Hall-matching refinement improves the
maximum-degree bound to $O(t^3(\log t)^{3/2})$ and, for connected graphs,
the order bound to $O(t^{19}(\log t)^9)$.
For the weaker local exclusion, a codegree graph and Ramsey's theorem give
the bound $\exp(O(t\log t))$. We also prove an exponential maximum-degree
bound when the full 1-subdivision $K_t^{(2)}$ is forbidden as a subgraph,
assuming minimum degree at least $t$. 
\end{abstract}

\section{Introduction}

In differential geometry, a manifold with positive Ricci curvature is
compact. Analogous results hold for graphs under suitable hypotheses. We
call a connected graph positively curved if its Lin--Lu--Yau curvature is
positive on every edge. Every connected positively curved graph with
bounded maximum degree is finite.

Lu and Wang proved that a planar graph with minimum degree at least three
and positive Lin--Lu--Yau curvature has maximum degree at most $17$, and
hence is finite \cite{LuWangPlanar}. For positively curved outerplanar
triangulations with minimum degree at least two, Brooks et al.
\cite{BrooksOuter} proved a maximum-degree bound of $9$. Liu, Lu, and Wang
\cite{LiuLuWangOuter} proved the sharp vertex bound $10$ for positively
curved outerplanar graphs.

Which topological restrictions, together with minimum-degree assumptions,
force positively curved graphs to have bounded maximum degree and hence be
finite? Planar graphs are precisely the graphs with no $K_5$ or $K_{3,3}$
minor, and outerplanar graphs are precisely the graphs with no $K_4$ or
$K_{2,3}$ minor. This motivates the study of positively curved
minor-free graphs.

Bai, Liu, and Song obtained
sharp maximum-degree bounds, together with order bounds, for positively
curved graphs with no $K_{2,t}$ minor \cite{BaiLiuSong}. Song recently
proved polynomial bounds for the maximum degree and order of positively
curved $K_t$-minor-free graphs \cite{SongDraft}:
\[
 \Delta(G)=O(t^5(\log t)^{3/2}),\qquad
 |V(G)|=O(t^{31}(\log t)^9).
\]

Song's proof uses edge-density estimates for minor-free graphs and
minor-preserving contractions. The matching argument below improves the
maximum-degree bound to $O(t^3(\log t)^{3/2})$ and, for connected graphs,
the order bound to $O(t^{19}(\log t)^9)$; see Theorem~\ref{thm:main}(c).

We use the following notation for topological minors. A subdivision of a
graph $H$ is obtained by replacing each edge of $H$ with a path so that the
replacement paths are internally vertex-disjoint. Write $\TK_t$ for the
family of all subdivisions of $K_t$, and write $\TKshort_t$ for the family
of subdivisions in which every replacement path has length at most two.
We call $G$ $\TK_t$-free if it contains no member of $\TK_t$, and
$\TKshort_t$-free if it contains no member of $\TKshort_t$. Thus ordinary
$\TK_t$-freeness is the stronger exclusion.
Write $K_t^{(2)}$ for the graph obtained by subdividing every edge of $K_t$
exactly once.

The following question was posed in \cite{SongDraft}:

\begin{problem}[\cite{SongDraft}, 4.1]
Let $t\ge5$. Is there a function $g(t)$ such that every connected graph
$G$ with no topological $K_t$ minor, minimum degree at least $t-1$, and
positive Lin--Lu--Yau curvature on every edge satisfies
$\Delta(G)\le g(t)$ and hence
$|V(G)|\le g(t)^{2g(t)(g(t)-1)}$?
\end{problem}
We answer this question positively and also prove a weaker bound for
$\TKshort_t$-free graphs.

\begin{theorem}\label{thm:main}
Let $t\ge 5$, and let $G$ be a locally finite graph with
$\delta(G)\ge t-1$ and positive Lin--Lu--Yau curvature on every edge.
\begin{enumerate}[label=\textup{(\alph*)}]
  \item If $G$ is $\TK_t$-free, then
  \[
       \Delta(G)\le C t^{7}.
  \]
  \item If $G$ is $\TKshort_t$-free, then
  \[
       \Delta(G)\le \exp\bigl(Ct\log t\bigr).
  \]
  \item If $G$ is $K_t$-minor-free, then
  \[
       \Delta(G)\le C t^3(\log t)^{3/2}.
  \]
  If $G$ is connected, then it also satisfies
  \[
       |V(G)|\le C t^{19}(\log t)^9.
  \]
\end{enumerate}
Here $C$ is an absolute constant. The exponents and constants are not
optimized.
\end{theorem}
The proof of part~(c) is given in the more general edge-density form of
Proposition~\ref{prop:minor-refinement} in Section~4. 

We also show that forbidding the single subgraph $K_t^{(2)}$, together
with the stronger minimum-degree assumption $\delta(G)\ge t$, suffices
for an exponential maximum-degree bound.
\begin{theorem}\label{thm:once-subdivision}
Let $t\ge5$, and let $G$ be a locally finite graph with
$\delta(G)\ge t$ and positive Lin--Lu--Yau curvature on every edge. If $G$
contains no copy of $K_t^{(2)}$ as a subgraph, then
\[
  \Delta(G)\le \exp\bigl(Ct\log t\bigr),
\]
where $C$ is an absolute constant.
\end{theorem}
The minimum-degree thresholds in Theorems~\ref{thm:main} and
\ref{thm:once-subdivision} are sharp. For $n\ge t$, the graph
$K_{t-2}\vee E_{n-t+2}$ is positively curved, $K_t$-minor-free, and has
minimum degree $t-2$ and maximum degree $n-1$ \cite{SongDraft}. The graph
$K_{t-1}\vee E_{n-t+1}$ is positively curved, $K_t^{(2)}$-free, and has
minimum degree $t-1$ and maximum degree $n-1$; its properties are proved
in Proposition~\ref{prop:once-threshold}.


The proof uses a common curvature-and-coloring argument. For the ordinary
topological exclusion, the theorem of Bollob\'as--Thomason and
Koml\'os--Szemer\'edi gives $O(t^2)$ edge density in every subgraph
\cite{BT,KS}, which
bounds a set of vertices with large codegree into one neighborhood. Under
the local exclusion this density statement is false, so a Ramsey argument
replaces it and yields an exponential bound.

Indeed, for $t\ge3$, girth at least seven rules out a $\TKshort_t$: the
paths corresponding to any three branch vertices would form a
cycle of length at most six. Erd\H{o}s' high-girth, high-chromatic
construction then gives $\TKshort_t$-free graphs with
arbitrarily large average degree, by passing to a subgraph of large minimum
degree \cite{Erdos}.

\begin{corollary}\label{cor:finite}
If $G$ satisfies the hypotheses of Theorem~\ref{thm:main} or
Theorem~\ref{thm:once-subdivision} and is connected, then $G$ is finite.
In each case its order satisfies
\[
 |V(G)|\le \Delta(G)^{2\Delta(G)(\Delta(G)-1)}.
\]
\end{corollary}

\begin{remark}
The local version is the stronger degree conclusion currently established
here only up to $\exp(O(t\log t))$. Obtaining a polynomial bound under
$\TKshort_t$-freeness is an open problem for this method.
\end{remark}

The rest of the paper is organized as follows. Section~2 recalls the
Lin--Lu--Yau curvature definition and the graph-density tools used later.
Section~3 develops the common curvature decomposition around a vertex.
Section~4 proves the ordinary topological-minor bound and its minor-free
refinement; Section~5 treats the locally topological exclusion, and
Section~6 treats exclusion of the full 1-subdivision $K_t^{(2)}$.
Section~7 compares the degree and order bounds and discusses open problems.
The proofs share a common strategy: curvature gives large neighborhoods in
an auxiliary graph, while density bounds or codegree estimates control the
exceptional sets. A clique-supported coloring then extracts a rainbow
clique, which lifts to the forbidden subdivision. In the minor-free case,
Hall's theorem matches the exceptional vertices to distinct pairs and
improves the degree estimate.

\section{Preliminaries}
All graphs are simple and locally finite. For $0\le\alpha<1$ and a vertex
$x$ with $d(x)>0$, let $m_x^\alpha$ be the probability measure that places mass
$\alpha$ at $x$ and distributes the remaining mass uniformly over $N(x)$:
\[
 m_x^\alpha(x)=\alpha,\qquad
 m_x^\alpha(z)=\frac{1-\alpha}{d(x)}\ (z\in N(x)),
\]
with zero mass elsewhere. For probability measures $\mu,\nu$ on $V(G)$,
their 1-Wasserstein distance is
\[
 W_1(\mu,\nu)=\min_{\pi\in\Pi(\mu,\nu)}
       \sum_{u,v\in V(G)} d(u,v)\pi(u,v),
\]
where $\Pi(\mu,\nu)$ is the set of couplings of $\mu$ and $\nu$.
For distinct vertices $x,y$ in the same component, define the Ollivier curvature
\[
 \kappa_\alpha(x,y)=1-\frac{W_1(m_x^\alpha,m_y^\alpha)}{d(x,y)}.
\]
The Lin--Lu--Yau curvature is the rescaled limit
\[
 \kappa(x,y)=\lim_{\alpha\uparrow1}
       \frac{\kappa_\alpha(x,y)}{1-\alpha};
\]
this limit exists for locally finite graphs \cite{LLY}. For an edge $xy$,
Kantorovich duality gives the useful formula
\[
 \kappa(x,y)=\inf_f\bigl(\Delta f(x)-\Delta f(y)\bigr),
 \qquad
 \Delta f(v)=\frac1{d(v)}\sum_{u\in N(v)}(f(u)-f(v)),
\]
where the infimum is over all functions $f:V(G)\to\mathbb R$ satisfying
$|f(u)-f(v)|\le1$ for every edge $uv$ and $f(y)-f(x)=1$. More generally,
for a geodesic
$x=v_0,v_1,\ldots,v_\ell=y$, the Wasserstein triangle inequality implies
\[
 \kappa(x,y)\ge\frac1\ell\sum_{i=1}^{\ell}
                    \kappa(v_{i-1},v_i).
\]
In particular, positive curvature on every edge implies positive curvature
between every pair of vertices in the same component. A connected graph $G$ is called
\emph{positively curved} if $\kappa(x,y)>0$ for every $xy\in E(G)$.
For $A\subseteq V(G)$ and
$v\in V(G)$, write
\[
 d_A(v)=|N(v)\cap A|,
 \qquad
 \codeg_A(u,v)=|N(u)\cap N(v)\cap A|.
\]
We use the following consequence of the topological clique theorem.
We write $R(a,b)$ for the usual two-color Ramsey number.

\begin{lemma}[Topological-clique density bound]\label{lem:density}
There is an absolute constant $C_0$ such that every finite
topological-$K_s$-free graph $H$ satisfies
\[
 e(H)\le C_0s^2|V(H)|.
\]
\end{lemma}

This follows from the theorem that average degree at least $C s^2$ forces a
subdivision of $K_s$ \cite{BT,KS}.

We use the following test-set inequality of Lu and Wang
\cite{LuWangPlanar}.

\begin{lemma}[Lu--Wang test-set inequality]\label{lem:LW}
Let $xy\in E(G)$, put $D=d(x)$, $d=d(y)$,
$C_{xy}=N(x)\cap N(y)$, and $c=|C_{xy}|$. For every nonempty
$S\subseteq N(y)\setminus\{x\}$, let $s=|S|$ and
$k=|S\cap N(x)|$. If $\kappa(x,y)>0$, then
\[
 \bigl|N(S)\cap(N(x)\setminus N(y))\bigr|
 >\frac{Ds}{d}-(k+c+1).
\]
\end{lemma}

\begin{lemma}[Clique-supported coloring]
\label{lem:supported-rainbow}
Let $F$ be a graph on $n\ge2t$ vertices with
$\delta(F)>(1-\rho)n$. For each edge $e$, give $e$ a color, and suppose
that every color $c$ has a supporting clique $Q_c$ with
$|Q_c|\le\sigma n$ such that every edge of color $c$ has both endpoints
in $Q_c$. Give each vertex $v$ a color, with all vertex colors distinct;
these colors may coincide with edge colors.
If
\begin{equation}\label{eq:rainbow}
 2\binom t2 \frac{1}{n} +
 2\binom t2\rho +3\binom t3 \sigma + 3\binom t3 \frac{1}{n} 
 <1,   
\end{equation}

then $F$ contains a rainbow $K_t$, i.e., all vertices and edges of this $K_t$ receive distinct colors.
\end{lemma}

\begin{proof}
Choose $t$ vertices independently and uniformly at random, allowing
repetitions. We consider four types of bad events.

\begin{enumerate}
    \item Some pair of selected vertices is equal. By a union bound, this
    has probability at most $\binom t2/n$.
    \item Some pair of distinct selected vertices is not an edge of $F$. This has
    probability at most $2\binom t2\rho$.
    \item There are distinct selected vertices $u,v,w$ with $uv\in E(F)$
    and $w\in Q_c$, where $c$ is the color of $uv$. This has probability
    at most $3\binom t3\sigma$.
    \item There are distinct selected vertices $u,v,w$ with $uv\in E(F)$
    and the color of $uv$ equal to the color of $w$. This has probability
    at most $3\binom t3/n$.
\end{enumerate}
If \eqref{eq:rainbow} holds, the union bound gives positive probability
that none of these bad events occurs. The first two conditions ensure that
the selected vertices induce a $K_t$. No two of its edges can have the same
color: if two did, choose one edge $uv$ and an endpoint $w$ of the other
edge with $w\notin\{u,v\}$. Then $w\in Q_c$, where $c$ is the color of
$uv$, so the third bad event occurs. The fourth condition prevents an edge
color from coinciding with a vertex color. Since the vertex colors are
already distinct, the resulting $K_t$ is rainbow.
\end{proof}

\section{A common curvature step}
The following neighborhood partition is due to Song \cite{SongDraft}. To
bound the degree $d(x)$ of each vertex $x$, it suffices to bound
$|U|$, $|W|$, and $|Z|$, defined as follows. 

Fix $x\in V(G)$ and write $D=d(x)$ and $A=N(x)$. Choose a parameter
$0<\eps<1$ and define
\[
 U=\{v\ne x:d_A(v)\ge\eps D\},
 \qquad
 W=\{y\in A\setminus U:|N(y)\cap U|\ge t-2\},
 \qquad
 Z=A\setminus(U\cup W).
\]
For $y\in Z$, put
\[
 R_y=N(y)\setminus(\{x\}\cup U),
 \qquad
 C_y=\{v\in A\setminus\{y\}:N(v)\cap R_y\ne\varnothing\}.
\]

\begin{figure}[htbp]
\centering
\begin{tikzpicture}[line cap=round,line join=round,every node/.style={font=\small}]
  \filldraw[fill=blue!4,draw=blue!65!black,line width=1pt]
    (0,0) ellipse (3.35 and 1.62);
  \filldraw[fill=green!12,draw=green!50!black,line width=0.8pt]
    (-1.55,0.48) ellipse (1.05 and 0.55);
  \filldraw[fill=violet!10,draw=violet!65!black,line width=0.8pt]
    (1.45,0.40) ellipse (1.35 and 0.8);
  \filldraw[fill=orange!15,draw=orange!75!black,line width=0.8pt]
    (-0.25,-1.33) ellipse (1.55 and 0.80);
  \draw[magenta!75!black,densely dotted,line width=1.1pt]
    (1.95,0.28) ellipse (0.73 and 0.42);

  \draw[gray!55,line width=0.45pt] (0,2.55)--(-2.05,0.36);
  \draw[gray!55,line width=0.45pt] (0,2.55)--(-1.55,0.62);
  \draw[gray!55,line width=0.45pt] (0,2.55)--(-0.95,0.22);
  \draw[gray!55,line width=0.45pt] (0,2.55)--(-0.35,-0.82);
  \draw[gray!55,line width=0.45pt] (0,2.55)--(0.82,0.88);
  \draw[gray!55,line width=0.45pt] (0,2.55)--(1.18,0.58);
  \draw[gray!55,line width=0.45pt] (0,2.55)--(1.90,0.34);
  \draw[gray!55,line width=0.45pt] (0,2.55)--(2.45,0.70);
  \draw[gray!55,line width=0.45pt] (0,2.55)--(-0.72,-1.47);
  \draw[gray!55,line width=0.45pt] (0,2.55)--(0.28,-1.52);
  \draw[gray!55,dashed] (4.05,-0.82) ellipse (0.38 and 0.30);

  \draw[violet!70!black,dashed,line width=0.9pt] (1.18,0.58)--(4.05,-0.82);
  \draw[violet!70!black,dashed,line width=0.9pt] (1.90,0.34)--(4.05,-0.82);

  \node[blue!65!black,anchor=west] at (-3.08,1.60) {$A=N(x)$};
  \node[green!45!black] at (-1.55,0.95) {$W$};
  \node[violet!65!black] at (1.43,1.02) {$Z$};
  \node[orange!75!black,anchor=east] at (-1.18,-1.45) {$U$};
  \node[magenta!75!black,anchor=west] at (2.25,-0.23) {$C_y$};
  \node at (0,2.78) {$x$};
  \node at (1.02,0.74) {$y$};
  \node[anchor=west] at (2.03,0.35) {$y'\in C_y$};
  \node[anchor=west] at (4.36,-0.89) {$z\in R_y$};

  \foreach \p in {(-2.05,0.36),(-1.55,0.62),(-0.95,0.22)}
    \fill[green!50!black] \p circle (1.6pt);
  \foreach \p in {(0.82,0.88),(1.18,0.58),(1.90,0.34),(2.45,0.70)}
    \fill[violet!65!black] \p circle (1.6pt);
  \foreach \p in {(-0.35,-0.82),(-0.72,-1.47),(0.28,-1.52),(0.25,-1.67)}
    \fill[orange!75!black] \p circle (1.6pt);
  \fill[gray!70!black] (4.05,-0.82) circle (1.6pt);
  \fill (0,2.55) circle (1.8pt);
\end{tikzpicture}
\caption{Neighborhood decomposition, following the schematic in Figure~1 of
\cite{SongDraft}. The set $A=N(x)$ is partitioned into
$(U\cap A)\mathbin{\dot\cup}W\mathbin{\dot\cup}Z$, while $U$ may also contain
vertices outside $A$. For the displayed $y\in Z$, the highlighted
$C_y\subseteq A\setminus\{y\}$ consists of vertices adjacent to $R_y$; the
dashed edges show one possible witness $z\in R_y$.}
\label{fig:neighborhood-sets}
\end{figure}
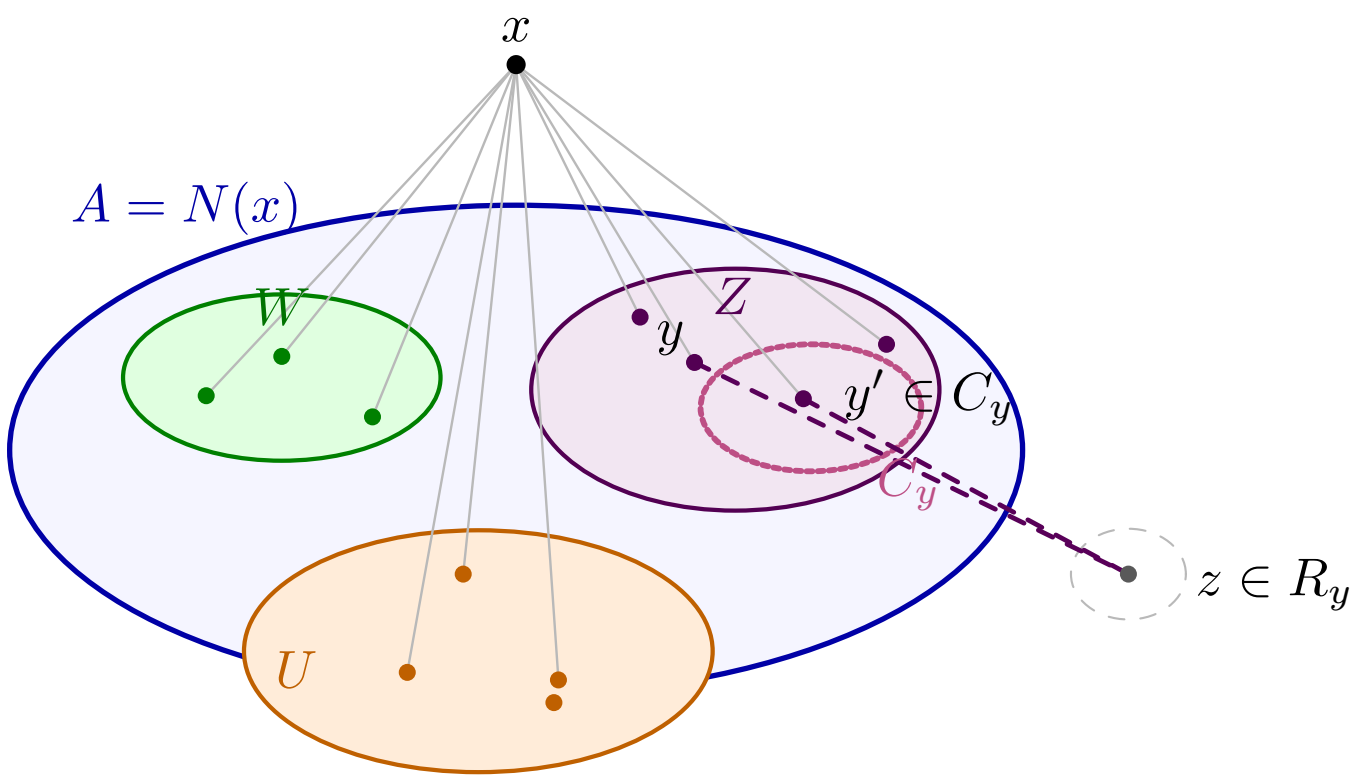

\begin{lemma}[Degree jump on $Z$]\label{lem:degree-jump}
Assume $\delta(G)\ge t-1$ and that every edge has positive curvature.
With $U,W,Z$ as above, for every $y\in Z$, writing
$c_y=|N(y)\cap A|$, we have $R_y\ne\varnothing$ and
\[
 d(y)>\frac{1}{2\eps+2/D}.
\]
\end{lemma}

\begin{proof}
Since $y\in Z$, we have $u_y:=|N(y)\cap U|\le t-3$. Thus
$|R_y|=d(y)-1-u_y\ge1$. Choose $z\in R_y$. Then $z\notin U$, so
$d_A(z)<\eps D$. Applying Lemma~\ref{lem:LW} to $xy$ with $S=\{z\}$,
and writing $k=|\{z\}\cap A|\le1$, gives
\[
 \eps D>|N(z)\cap(A\setminus N(y))|
 >\frac{D}{d(y)}-(k+c_y+1).
\]
Therefore $D/d(y)<\eps D+k+c_y+1\le\eps D+c_y+2<2\eps D+2$.
Rearranging yields the degree lower bound.
\end{proof}

\begin{lemma}\label{lem:curv-neighborhood}
Assume $\delta(G)\ge t-1$, all edges have positive curvature, and
$U,W,Z$ are as above. Then $R_y\ne\varnothing$ for each $y\in Z$, and
\[
 |C_y|>(1-2(t-1)\eps)D-2t+2.
\]
\end{lemma}

\begin{proof}
Let $c_y=|N(y)\cap A|$. Since $y\notin U$, $c_y<\eps D$. Also
$u_y:=|N(y)\cap U|\le t-3$. By Lemma~\ref{lem:degree-jump},
$R_y\ne\varnothing$ and $D/d(y)<2\eps D+2$.
Now apply Lemma~\ref{lem:LW} with $S=R_y$. Since
$d(y)-|R_y|=1+u_y\le t-2$ and $|R_y\cap A|\le c_y$,
\begin{align*}
 |N(R_y)\cap(A\setminus N(y))|
 &> \frac{D|R_y|}{d(y)}-(|R_y\cap A|+c_y+1)\\
 &\ge D-(t-2)\frac D{d(y)}-2c_y-1\\
 &>D-(t-2)(2\eps D+2)-2\eps D-1\\
 &=(1-2(t-1)\eps)D-2t+3.
\end{align*}
Since $t-2>0$, substituting the upper bound for $D/d(y)$ into the term
with coefficient $-(t-2)$ reverses the inequality. The additive
$2(t-2)$ is retained in the constant term.
The set on the left contains $y$; all its other vertices belong to $C_y$.
Removing $y$ proves the claim.
\end{proof}

\begin{lemma}\label{lem:colored-Z}
Assume $\delta(G)\ge t-1$ and that every edge has positive curvature.
Use the sets $A,U,W,Z$ defined above with
$\epsilon=\frac{1}{Ct^3}$, where $C\geq 32$ is a constant, and put
$B=|A\setminus Z|$. If
\[
 D\geq 2B+32t^3,
\]
then $G$ contains a subgraph $K^{(2)}_t$.
\end{lemma}

\begin{proof}
Let $F$ have vertex set $Z$. Join $y,y'\in Z$ if they have a common neighbor in $R_y\cap R_{y'}$.
Give each vertex a color of itself.
For every edge, choose one common neighbor
as its color. 
For a witness color $z$ of edges, all edges of color $z$ lie in the
clique $Q_z=N(z)\cap Z$ of $F$. Since $z\notin U$,
$|Q_z|\le |N(z)\cap A|<\eps D$. 

Every $y\in Z$ has at least $|C_y|-B$ neighbors in $F$. Lemma
\ref{lem:curv-neighborhood} and $|Z|=D-B$ imply
\[
 \delta(F)>(1-\rho)|Z|,
 \qquad
\rho\le 2\bigl(2(t-1)\eps+(2t-2)/D\bigr),
\]
since $D\ge2B$. Also $n:=|Z|=D-B\ge D/2$.
The witness supports have size at most $2\eps n$.  The defect bound in Lemma
\ref{lem:supported-rainbow} is less than one, with $\sigma=2\eps$.
Indeed, its terms are
\[
2\binom t2 \frac{1}{n} +
 2\binom t2\rho +3\binom t3 \sigma +3\binom t3 \frac{1}{n} 
 <16t^3\left(\epsilon +\frac{1}{D}\right).
\]

Lemma~\ref{lem:supported-rainbow} gives a rainbow $K_t$ in $F$ whose
witness colors are not branch vertices. Replace each witness-colored
edge by its length-two path through its witness, giving a subgraph $K^{(2)}_t$ in $G$.
\end{proof}

\section{Ordinary topological exclusion}

\begin{lemma}[Hall matching for $W$]\label{lem:W-pair-graph}
Fix $x$ and the associated sets $A,U,W$. Suppose $G$ is locally finite and
$\TK_t$-free, and
\[
 \eps D\ge \binom{t-2}{2}+2t.
\]
Then there is a $\TK_t$-free graph $H$ on $U$ with
$e(H)=|W|$. If $G$ is $K_t$-minor-free, then $H$ can be taken
$K_t$-minor-free as well.
\end{lemma}

\begin{proof}
Put $q=\binom{t-1}{2}+1$. For each $y\in W$, choose a set
$T_y\subseteq N(y)\cap U$ of size $t-2$, and let $P_y$ consist of the
pairs $uv\in\binom{T_y}{2}$ such that $uv\notin E(G)$ and
$\codeg_{A\setminus U}(u,v)<q$. This set is nonempty: otherwise take
$x,y,T_y$ as branch vertices. Use direct edges for adjacent pairs in
$T_y$, and route the remaining pairs through distinct common neighbors in
$A\setminus U$. At each common-neighbor choice at most
$1+(t-2)+\binom{t-2}{2}-1=q-1$ vertices are forbidden. Then route each
$xu$ through a distinct vertex of $A\cap N(u)$, avoiding all branch
vertices and all pair-path interiors. At most
$\binom{t-2}{2}+2t-4$ vertices are forbidden at any such choice, so the
degree assumption supplies the required interiors. This gives a subdivision
of $K_t$, a contradiction.

Since $W\subseteq A$ and $G$ is locally finite, $W$ is finite.
Form a bipartite graph with left side $W$ and right side the pairs in
$\bigcup_{y\in W}P_y$, joining $y$ to the pairs in $P_y$. Suppose Hall's
condition fails, and choose an inclusion-minimal $S\subseteq W$ with
$|N(S)|<|S|$. Write $h=|N(S)|$; minimality gives $h=|S|-1$. Fix
$y\in S$. For $T\subseteq P_y$, let $S_T$ be the vertices of $S$ joined
to at least one pair in $T$. If $|S_T\setminus\{y\}|<|T|$, then
$y\in S_T$, because $y$ is joined to every pair in $P_y$, and hence
$|S_T|\le|T|$. Since
$|T|\le|N(S)|=|S|-1$, the set $S'=S\setminus S_T$ is nonempty. Also
$N(S')\subseteq N(S)\setminus T$, so
\[
 |N(S')|\le h-|T|<|S|-|T|\le |S'|,
\]
contrary to the minimality of $S$. Hence Hall's condition holds between
$P_y$ and $S\setminus\{y\}$, so the pairs in $P_y$ can be matched to
distinct vertices of $S\setminus\{y\}$. Write this matching as
$\phi:P_y\to S\setminus\{y\}$. For every $e=uv\in P_y$, we have
$e\in P_{\phi(e)}$, so $\phi(e)\in A\setminus U$ is a common neighbor
of $u$ and $v$. These distinct vertices serve as the interiors for the
pairs in $P_y$.

Take $x,y,T_y$ as branch vertices. Route every pair $e\in P_y$ through
$\phi(e)$. For every remaining nonadjacent pair
$uv\in\binom{T_y}{2}\setminus P_y$, choose a common neighbor in
$A\setminus U$. At each choice at most
\[
 1+|P_y|+(t-2)+\left(\binom{t-2}{2}-|P_y|-1\right)
 =q-1
\]
vertices are forbidden, so the codegree threshold supplies a distinct
interior. Now, for each $u\in T_y$, choose a distinct interior for the
$xu$ path in $A\cap N(u)$, avoiding all branch vertices and all pair-path
interiors. At most
$\binom{t-2}{2}+2t-4$ vertices are forbidden, so the assumption on
$\eps D$ guarantees a choice. Along with $xy$ and the edges from $y$ to
$T_y$, these paths form a subdivision of $K_t$, a contradiction. Thus no
inclusion-minimal Hall-deficient set exists. Since $W$ is finite,
Hall's condition holds globally, and there is a matching of all vertices
of $W$ to distinct pairs in $U$.

Let $H$ have vertex set $U$ and one edge for each matched pair, labeled
by its matched vertex $y\in W$. Let $H^*$ be obtained by subdividing each
edge $uv$ once using its label $y$. The labels are distinct, so $H^*$ is a
subgraph of $G$. Any subdivision of $K_t$ in $H$ lifts to one in $H^*$,
so $H$ is $\TK_t$-free. If $G$ is $K_t$-minor-free, then $H$ is also
$K_t$-minor-free, since contracting the subdivided edges in $H^*$ recovers
$H$. In either case $e(H)=|W|$.
\end{proof}

\begin{proof}[Proof of Theorem~\ref{thm:main}(a)]
Let $a_t=C_0t^2$, where $C_0$ is as in Lemma~\ref{lem:density}, and set
$\eps=1/(Ct^3)$ with $C$ sufficiently large. If
$D<\max\{4a_t/\eps, (\binom{t-2}{2}+2t)/\eps\}$, then $D=O(t^5)$ and
there is nothing to prove.

Otherwise, for every finite $U_0\subseteq U$, apply
Lemma~\ref{lem:density} to $G[A\cup U_0]$:
\[
 \eps D|U_0|
 \le \sum_{u\in U_0}d_A(u)
 \le 2e(G[A\cup U_0])
 \le 2a_t(D+|U_0|).
\]
Since $D\ge4a_t/\eps$, this gives $|U_0|\le4a_t/\eps$, and hence
$|U|=O(t^5)$. The bound for all finite $U_0$ also proves that $U$ is
finite.

Apply Lemma~\ref{lem:W-pair-graph}. The graph $H$ is
 $\TK_t$-free, so Lemma~\ref{lem:density} yields
\[
 |W|=e(H)\le C_0t^2|U|=O(t^7).
\]
Consequently $B=|A\setminus Z|\le |U|+|W|=O(t^7)$.
If $D<2B+32t^3$, then $D=O(t^7)$ and we are done.
Lemma~\ref{lem:colored-Z} contradicts $K^{(2)}_t$-freeness, which is
implied by ordinary $\TK_t$-freeness. Thus
$D=O(t^7)$ uniformly over all vertices $x$.
\end{proof}

\begin{proposition}[Minor-free refinement]\label{prop:minor-refinement}
Let $t\ge5$, and let $G$ be a locally finite, positively curved graph
with $\delta(G)\ge t-1$ and no $K_t$ minor. Suppose every finite
$K_t$-minor-free graph $H$ satisfies $e(H)\le c|V(H)|$, where
$c\ge t-2$. Then
\[
  \Delta(G)=O(c^3),
\]
with an absolute implied constant. In particular, the Kostochka--Thomason
theorem \cite{Kostochka,ThomasonMinor} gives $c=O(t\sqrt{\log t})$ and
hence $\Delta(G)=O(t^3(\log t)^{3/2})$. If $G$ is connected, Song's order
estimate also gives $|V(G)|=O(t^{19}(\log t)^9)$ \cite{SongDraft}.
\end{proposition}

\begin{proof}
Fix $x\in V(G)$, and use the common neighborhood decomposition with
$A=N(x)$, $D=d(x)$, and $\eps=1/(5c)$. It suffices to consider
$D\ge25c^3$. Put $e=e(G[A])$, $U_1=U\cap A$, and $U_2=U\setminus A$.
For every finite $T\subseteq V(G)\setminus N[x]$, applying the edge bound
to $G[A\cup T\cup\{x\}]$ gives
\[
 \sum_{z\in T}(d_A(z)-c)\le c(D+1)-e-D.
\]
Applying the edge bound to $G[U_1]$ yields
\[
 e\ge \sum_{u\in U_1}d_A(u)-e(G[U_1])
   \ge (\eps D-c)|U_1|.
\]
For every finite $T\subseteq U_2$, we have $d_A(z)\ge\eps D$ for
$z\in T$, and therefore
\[
 (\eps D-c)|T|\le c(D+1)-e-D.
\]
This uniform bound shows that $U_2$ is finite and that the same inequality
holds with $T=U_2$. Combining it with the preceding estimate for $e$ gives
\[
 (\eps D-c)|U|\le(c-1)D+c.
\]
Since $D\ge25c^3$, it follows that
\[
 |U|\le \frac{5c((c-1)D+c)}{D-5c^2}\le5c^2-4c.
\]

Moreover, $\eps D\ge5c^2\ge\binom{t-2}{2}+2t$, so
Lemma~\ref{lem:W-pair-graph} applies. Its lifted graph $H$ is
$K_t$-minor-free, and hence
\[
 |W|=e(H)\le c|U|.
\]

Song's Claims 3 and 4 in the proof of Theorem~1.3 of \cite{SongDraft} give
the following estimates, where $c_y=|N(y)\cap A|$:
\[
 |C_y|>\frac{4D}{5}-t(c_y+3),
 \qquad
 \sum_{y\in Z}|C_y|\le\frac{2D^2}{5}.
\]
The auxiliary graphs in Claim 4 are $K_{t-1}$-minor-free, hence
$K_t$-minor-free, so the assumed edge bound applies to them with this same
$c$.
In Claim 4, the lower bound on $D$ is used to ensure
$r=\lfloor\eps D\rfloor\ge2$, which holds here. Since
$\sum_{y\in Z}c_y\le2e(G[A])\le2cD$, we obtain
\[
 |Z|\left(\frac{4D}{5}-3t\right)
 \le\frac{2D^2}{5}+2tcD.
\]
Here $D\ge25c^3$ and $c\ge t-2$ imply $D\ge15t$. Hence
\[
 |Z|\le\frac D2+
 \frac{(2tc+3t/2)D}{4D/5-3t}
 \le\frac D2+4tc+3t.
\]
Finally, $D=|A|\le|U|+|W|+|Z|$, so
\[
 D\le(1+c)(5c^2-4c)+\frac D2+4tc+3t,
\]
which gives $D=O(c^3)$. If $D<25c^3$, the same conclusion is immediate.
The bound is uniform in $x$. By the Kostochka--Thomason theorem, every
finite $K_t$-minor-free graph has $e(H)\le c|V(H)|$ with
$c=O(t\sqrt{\log t})$. This proves Theorem~\ref{thm:main}(c).

If $G$ is connected, Song's estimate
$|V(G)|\le 2t\Delta(G)^6$ (Theorem~1.3 in \cite{SongDraft}) yields the
stated order bound.

\end{proof}

\section{Locally topological exclusion}

The set $U$ need not have bounded size from edge density in this case:
$\TKshort_t$-free graphs can have unbounded average degree.
Instead we use codegrees to build a graph on $U$.

\begin{lemma}\label{lem:local-U}
Let $G$ be $K_t^{(2)}$-free, and assume
$\eps^2D>2Q$, where $Q=\binom{t+1}{2}$. Let $J$ be the
graph on $U$ in which $uv\in E(J)$ if and only if $\codeg_A(u,v)\ge Q$.
Then $J$ is $K_{t-1}$-free and
\[
 |U|<R\left(t-1,\left\lceil\frac2\eps\right\rceil\right)
 \le \binom{t+\lceil2/\eps\rceil-3}{t-2}.
\]
\end{lemma}

\begin{proof}
Work first with an arbitrary finite subset $U_0\subseteq U$, and let $J_0$
be the graph on $U_0$ in which $uv\in E(J_0)$ precisely when
$\codeg_A(u,v)\ge Q$. We first show that $J_0$ is $K_{t-1}$-free. If
$u_1,\dots,u_{t-1}$ formed a clique in $J_0$, take
$x,u_1,\dots,u_{t-1}$ as
branch vertices. For each $i$, choose a distinct common neighbor
$w_i\in A\cap N(u_i)$ for the pair $xu_i$, avoiding the branch vertices.
This is possible greedily because $|A\cap N(u_i)|\ge\eps D\ge2t$.
For each pair $u_i,u_j$, choose a distinct common neighbor in $A$ for the
pair path. At most $t-1$ branch vertices, $t-1$ previously chosen
$x$-path interiors, and $\binom{t-1}{2}-1$ earlier pair interiors are
excluded. This is fewer than $Q$, so the codegree threshold supplies a
choice at every step. These paths form a copy of $K_t^{(2)}$, a
contradiction.

Let $J_0^c$ be the complement of $J_0$. Its independence number is at most
$t-2$. If $X$ is a clique of order $h$ in $J_0^c$, then the sets
$S_u=N(u)\cap A$, $u\in X$, have sizes at least $\eps D$ and pairwise
intersections at most $Q-1$. Put $s=\sum_{u\in X}|S_u|$ and
$d_a=|\{u\in X:a\in S_u\}|$. Then
\[
 \sum_{a\in A}\binom{d_a}{2}
 \ge \frac12\left(\frac{s^2}{D}-s\right),
 \qquad
 \sum_{a\in A}\binom{d_a}{2}
 \le (Q-1)\binom h2.
\]
If $h\ge2/\eps$, these inequalities imply $\eps^2D\le2Q$, contrary to
the hypothesis. Indeed, $s\ge h\eps D\ge2D$, so the lower bound is at
least $h^2\eps^2D/4$, while the upper bound is less than $Qh^2/2$.
Hence $\omega(J_0^c)<\lceil2/\eps\rceil$. Ramsey's theorem
now gives the stated bound. Since the same bound holds for every finite
$U_0\subseteq U$, it also implies that $U$ itself is finite and satisfies
this bound. Every $K_{t-1}$ in $J$ would lie in some such $J_0$, so $J$ is
$K_{t-1}$-free as well.
\end{proof}

\begin{proof}[Proof of Theorem~\ref{thm:main}(b)]
Set $Q=\binom{t+1}{2}$ and $\eps=1/(Kt^3)$, with $K$ sufficiently large.
If $D\le\max\{4Q/\eps^2,2t/\eps\}$, then $D=O(t^{8})$. Otherwise,
Lemma~\ref{lem:local-U} applies, and
\[
 |U|\le R(t-1,\lceil2/\eps\rceil)=\exp(O(t\log t)).
\]

As in the ordinary case, if a vertex $y\in W$ had $t-2$ neighbors in $U$
all pairwise adjacent in $J$, the codegree threshold $Q$ would provide
distinct interiors for the pair paths and give a $\TKshort_t$. Thus every
$y\in W$ is chargeable to a pair not in $J$, and
\[
 |W|\le (Q-1)\binom{|U|}{2}.
\]
It follows that $B=|A\setminus Z|\le |U|+|W|
=\exp(O(t\log t))$.

If $D<2B+32t^3$, then $D=\exp(O(t\log t))$ and we are done.
Otherwise, Lemma~\ref{lem:colored-Z} gives a subgraph $K^{(2)}_t$, again a
contradiction. Combining the cases proves
$D\le\exp(O(t\log t))$ for every vertex.
\end{proof}

\section{Excluding the full 1-subdivision}

We now forbid only the graph $K_t^{(2)}$, in which every edge of $K_t$ is
subdivided exactly once. This exclusion allows direct edges among the
branch vertices, so the dense graph used below records only pairs with a
common-neighbor witness. The minimum-degree assumption is raised to $t$.

Fix $x\in V(G)$, let $D=d(x)$ and $A=N(x)$, and choose $\eps>0$. Define
\[
 U=\{v\ne x:d_A(v)\ge\eps D\},\qquad
 W_2=\{y\in A\setminus U:|N(y)\cap U|\ge t-1\},\qquad
 Z_2=A\setminus(U\cup W_2).
\]
For $y\in Z_2$, put
\[
 R_y^{(2)}=N(y)\setminus(\{x\}\cup U),\qquad
 C_y^{(2)}=\{v\in A\setminus\{y\}:N(v)\cap R_y^{(2)}\ne\varnothing\}.
\]

\begin{lemma}[Neighborhood expansion for $Z_2$]\label{lem:once-neighborhood}
Assume $\delta(G)\ge t$ and positive curvature on every edge. Then for
every $y\in Z_2$, $R_y^{(2)}\ne\varnothing$ and
\[
 |C_y^{(2)}|>(1-2t\eps)D-2t.
\]
\end{lemma}

The proof follows the same argument as Lemma~\ref{lem:curv-neighborhood},
using $|N(y)\cap U|\le t-2$ and $d(y)\ge t$. We omit the details.

\begin{lemma}[Coloring lemma for $Z_2$]\label{lem:colored-Z2}
Assume $\delta(G)\ge t$ and that every edge has positive curvature. Use
the sets $A,U,W_2,Z_2$ above with $\eps=1/(Ct^3)$, where $C\ge32$, and put
$B=|A\setminus Z_2|$. If
\[
 D\ge2B+32t^3,
\]
then $G$ contains a subgraph $K_t^{(2)}$.
\end{lemma}

The proof is the same as that of Lemma~\ref{lem:colored-Z}, using
$R_y^{(2)}$, $C_y^{(2)}$, and Lemma~\ref{lem:once-neighborhood}; we omit
it.

\begin{proof}[Proof of Theorem~\ref{thm:once-subdivision}]
Fix $x\in V(G)$ and use the notation above. Set
$Q=\binom{t+1}{2}$ and $\eps=1/(Ct^3)$. If
$D\le\max\{2t/\eps,2Q/\eps^2\}$, then $D=O(t^{8})$. Otherwise
Lemma~\ref{lem:local-U} applies, and
\[
 |U|<R(t-1,\lceil2/\eps\rceil)
      =\exp(O(t\log t)).
\]

Let $J$ be the codegree graph on $U$ from Lemma~\ref{lem:local-U}.
For each $y\in W_2$, choose $t-1$ neighbors in $U$. They cannot form a
clique in $J$, since $J$ is $K_{t-1}$-free. Thus some pair $u,u'\in U$ has
$\codeg_A(u,u')<Q$. Assign $y$ to one such pair. For a fixed pair,
at most $Q-1$ vertices of $A$ are common neighbors of $u,u'$, and hence
\[
 |W_2|\le(Q-1)\binom{|U|}{2}.
\]
Consequently $B=|A\setminus Z_2|\le |U|+|W_2|
=\exp(O(t\log t))$.

If $D<2B+32t^3$, then $D=\exp(O(t\log t))$ and we are done.
Otherwise, Lemma~\ref{lem:colored-Z2}
produces a copy of $K_t^{(2)}$, a contradiction. The bound holds for
every choice of $x$, proving the theorem.
\end{proof}

For $a,b\ge1$, write $K^*_{a,b}=K_a\vee E_b$, where $E_b$ is the
edgeless graph on $b$ vertices.
\begin{proposition}[Sharpness of the minimum-degree hypothesis]
\label{prop:once-threshold}
For every $t\ge5$ and $n\ge t$, the graph
$K^*_{t-1,n-t+1}=K_{t-1}\vee E_{n-t+1}$ is positively curved,
$K_t^{(2)}$-free, and satisfies
\[
 \delta(K^*_{t-1,n-t+1})=t-1,
 \qquad
 \Delta(K^*_{t-1,n-t+1})=n-1.
\]
Consequently, the lower-degree hypothesis $\delta(G)\ge t$ in
Theorem~\ref{thm:once-subdivision} cannot be weakened to $t-1$.
\end{proposition}

\begin{proof}
Let $S$ be the $(t-1)$-vertex clique and $I$ the $(n-t+1)$-vertex
independent set. Each vertex of $S$ has degree $n-1$, and each vertex of
$I$ has degree $t-1$. The two edge types have Lin--Lu--Yau curvatures
\[
 \kappa(u,v)=\frac{n}{n-1}\quad(u,v\in S),
 \qquad
 \kappa(u,w)=\frac{t}{n-1}\quad(u\in S,\ w\in I),
\]
as follows from the LLY curvature formula (see also
Proposition~1.2 of \cite{SongDraft}). Thus all edges have
positive curvature.

It remains to show that the graph contains no copy of $K_t^{(2)}$. Suppose
there were such a copy, and let $B$ be its set of $t$ branch vertices.
Since $|S|=t-1$, at least one branch vertex lies in $I$. If
$1\le |B\cap I|\le t-1$, choose $w\in B\cap I$. The $t-1$ distinct
subdivision vertices incident with $w$ must all lie in $S\setminus B$,
but $|S\setminus B|=|B\cap I|-1\le t-2$, a contradiction. If
$B\subseteq I$, then all $\binom{t}{2}$ subdivision vertices must lie
in $S$, which is impossible since $\binom{t}{2}>t-1=|S|$.
Therefore the graph is $K_t^{(2)}$-free.
\end{proof}

\section{Order bounds and open directions}

For a connected graph of positive Lin--Lu--Yau curvature and maximum
degree $\Delta\ge3$, the standard lower bound
$\kappa(x,y)\ge1/(\Delta(\Delta-1))$ for every edge $xy$, together with
the discrete Bonnet--Myers theorem, gives
\[
 \operatorname{diam}(G)\le2\Delta(\Delta-1),
 \qquad
 |V(G)|\le\Delta^{2\Delta(\Delta-1)}.
\]
This proves Corollary~\ref{cor:finite}.

We compare the degree and order bounds for the four exclusion conditions
considered in this paper and in Song's minor-free result. We assume
$t\ge5$, connectedness, positive Lin--Lu--Yau curvature on every edge, and
the minimum degree shown in the table.

\begin{table}[htbp]
\centering
\small
\setlength{\tabcolsep}{3pt}
\renewcommand{\arraystretch}{1.15}
\begin{tabular}{|p{2cm}|p{3cm}|p{3cm}|p{3cm}|p{2.7cm}|}
\hline
 & \shortstack{$K_t^{(2)}$-subgraph-free}
 & \shortstack{$\TKshort_t$-free}
 & \shortstack{$\TK_t$-free}
 & \shortstack{$K_t$-minor-free} \\
\hline
Minimum degree
 & $\delta\ge t$
 & $\delta\ge t-1$
 & $\delta\ge t-1$
 & $\delta\ge t-1$ \\
\hline
Maximum degree
 & $\exp(O(t\log t))$
 & $\exp(O(t\log t))$
 & $O(t^7)$
 & $O(t^3(\log t)^{3/2})$ (Theorem~\ref{thm:main}(c)) \\
\hline
Vertex bound
 & $\exp(\exp(O(t\log t)))$
 & $\exp(\exp(O(t\log t)))$
 & $\exp(O(t^{14}\log t))$
 & $O(t^{19}(\log t)^9)$ (Theorem~\ref{thm:main}(c)) \\
\hline
\end{tabular}
\caption{Degree and order bounds under four clique-exclusion conditions.}
\label{tab:exclusion-bounds}
\end{table}

For the first three cases, the order bounds follow by substituting the
maximum-degree estimates into the curvature volume bound
$|V(G)|\le\Delta^{2\Delta(\Delta-1)}$. In particular, the proof for
$\TK_t$-free graphs gives a polynomial bound on $\Delta$ but only the
exponential order bound shown in the table; it does not establish a
polynomial order bound in $t$. The $O(t^2)$ edge-density theorem for
topological-clique-free graphs controls the sets used in the local
curvature argument and yields the polynomial degree estimate. It does not
provide the additional global structure needed for a polynomial order
estimate.

The minor-free case has stronger global tools. Song's minor-specific
contraction and treewidth arguments give the direct order estimate
$|V(G)|<2t\Delta(G)^6$ \cite{SongDraft}. Theorem~\ref{thm:main}(c)
improves the degree bound to $O(t^3(\log t)^{3/2})$; substituting it into
Song's order estimate gives $|V(G)|=O(t^{19}(\log t)^9)$. 
For the two
Ramsey-based exclusions, the present argument gives an exponential degree
bound and consequently a double-exponential order bound. The table records
the bounds established here;
in particular, the non-polynomial order bound in the topological case is
a limitation of the current argument, not a lower bound ruling out a
polynomial one.

The principal unresolved issue is polynomial dependence under the weaker
exclusions. Ordinary topological exclusion supplies linear edge bounds in
every subgraph, while local topological and $K_t^{(2)}$ exclusion do not.
The Ramsey bound on $|U|$ is the source of the exponential loss in both
arguments. A polynomial estimate for $|U|$, or a direct curvature bound on
$|U|$, would lead to polynomial maximum-degree bounds through the same
charging and colored-neighborhood steps.

\begin{remark}
For $t=3,4$, having no topological $K_t$ minor is equivalent to having
no $K_t$ minor. For $t=3$, $\TKshort_3$-freeness is equivalent to having
no cycle of length at most six. The only connected positively curved
graphs with this property are stars, so none satisfies $\delta(G)\ge2$.
The remaining small cases for $t=3,4$ require separate treatment.
\end{remark}

\section*{Acknowledgements}
The author thanks Dr.~Zi-Xia Song for communicating the problem and ChatGPT-6 (Luna) for assistance with proofreading and exploring proof strategies during the preparation of this manuscript. The author takes full responsibility for the contents of the paper.


\begin{thebibliography}{9}

\bibitem{BT}
B. Bollob\'as and A. Thomason,
\newblock Proof of a conjecture of Mader, Erd\H{o}s and Hajnal on topological
complete subgraphs,
\newblock \emph{European Journal of Combinatorics} 19 (1998), 883--887.

\bibitem{KS}
J. Koml\'os and E. Szemer\'edi,
\newblock Topological cliques in graphs II,
\newblock \emph{Combinatorics, Probability and Computing} 5 (1996), no. 1,
79--90.

\bibitem{Kostochka}
A. V. Kostochka,
\newblock Lower bound of the Hadwiger number of graphs by their average
degree,
\newblock \emph{Combinatorica} 4 (1984), no. 4, 307--316.

\bibitem{ThomasonMinor}
A. Thomason,
\newblock An extremal function for contractions of graphs,
\newblock \emph{Mathematical Proceedings of the Cambridge Philosophical
Society} 95 (1984), no. 2, 261--265.

\bibitem{LLY}
Y. Lin, L. Lu, and S.-T. Yau,
\newblock Ricci curvature of graphs,
\newblock \emph{Tohoku Mathematical Journal} 63 (2011), 605--627.

\bibitem{LuWangPlanar}
L. Lu and Z. Wang,
\newblock On the size of planar graphs with positive Lin--Lu--Yau Ricci
curvature,
\newblock \emph{European Journal of Combinatorics} 135 (2026), article
104355, arXiv:2010.03716.

\bibitem{BrooksOuter}
G. Brooks, F. Osaye, A. Schenfisch, Z. Wang, and J. Yu,
\newblock Outerplanar graphs with positive Lin--Lu--Yau curvature,
\newblock \emph{Journal of Combinatorics} 16 (2025), no. 4, 465--480,
\newblock doi:10.4310/JOC.250923000902, arXiv:2403.04110.

\bibitem{LiuLuWangOuter}
X. Liu, L. Lu, and Z. Wang,
\newblock On the size of outerplanar graphs with positive Lin--Lu--Yau Ricci
curvature,
\newblock arXiv:2409.13666, 2024.

\bibitem{SongDraft}
Z.-X. Song,
\newblock On the maximum degree and order of $K_t$-minor-free graphs with
positive Lin--Lu--Yau curvature,
\newblock arXiv:2610.08706, version 3, 2026.

\bibitem{BaiLiuSong}
S. Bai, X. Liu, and Z.-X. Song,
\newblock Positively Lin--Lu--Yau curved graphs with no $K_{2,t}$ minor,
\newblock arXiv:2610.05519, 2026.

\bibitem{Erdos}
P. Erd\H{o}s,
\newblock Graph theory and probability,
\newblock \emph{Canadian Journal of Mathematics} 11 (1959), 34--38.

\end{thebibliography}
\end{document}